\documentclass[11pt,a4paper,leqno]{amsart}
\usepackage[utf8]{inputenc}
\usepackage[T1]{fontenc}
\usepackage[english]{babel}
\usepackage{lmodern}
\usepackage{amsmath,amssymb,amsthm,mathtools}
\usepackage[a4paper,left=29mm,right=29mm,top=27mm,bottom=28mm]{geometry}
\usepackage{microtype}
\usepackage{titlesec}
\usepackage{hyperref}
\hypersetup{hidelinks,unicode=true,
 pdftitle={Density of polynomials whose Aron--Berner extensions attain their norms},
 pdfsubject={A polynomial Lindenstrauss theorem over the real and complex fields},
 pdfauthor={Javier Falc\'o}}

\newtheorem{theorem}{Theorem}[section]
\newtheorem{lemma}[theorem]{Lemma}

\theoremstyle{remark}
\newtheorem{remark}[theorem]{Remark}
\theoremstyle{definition}
\titleformat{\section}[block]{\centering\normalfont\bfseries}{\S\thesection.}{.65em}{}
\titlespacing*{\section}{0pt}{3.2ex plus 1ex minus .2ex}{1.7ex plus .3ex}
\makeatletter
\def\@settitle{\begin{center}\bfseries\Large\@title\end{center}}
\makeatother
\allowdisplaybreaks[2]
\newcommand{\K}{\mathbb K}
\newcommand{\R}{\mathbb R}
\newcommand{\C}{\mathbb C}
\newcommand{\N}{\mathbb N}
\newcommand{\T}{\mathbb T}
\newcommand{\Pn}{\mathcal P({}^nX)}
\newcommand{\PnY}{\mathcal P({}^nX;Y)}

\newcommand{\Rea}{\operatorname{Re}}
\newcommand{\conv}{\operatorname{conv}}
\newcommand{\Rop}{\mathcal R}
\newcommand{\norm}[1]{\left\lVert#1\right\rVert}

\newcommand{\eps}{\varepsilon}
\newcommand{\U}{\mathcal U}
\newcommand{\ul}[1]{\lim_{#1\to\U}}

\title[Norm attaining Aron--Berner extensions]{Density of polynomials whose Aron--Berner extensions attain their norms}
\author{Javier Falc\'o}
\address[Javier Falc\'o]{Departamento de An\'alisis Matem\'atico, Universidad de Valencia, Doctor Moliner 50, 46100 Burjasot (Valencia), Spain}
\email{francisco.j.falco@uv.es}
\thanks{Supported by Project PID2025-168101NB-I00 funded by MICIU/AEI/10.13039/501100011033.}
\date{}

\dedicatory{Dedicated to Professors Domingo García and Manuel Maestre on the occasion of their 70th birthdays.}

\begin{document}
\begin{abstract}
We establish a polynomial analogue of Lindenstrauss' theorem in full generality. More precisely, we prove that, for every real or complex Banach space \(X\) and every positive integer \(n\), the \(n\)-homogeneous scalar-valued polynomials whose Aron–Berner extensions to \(X^{**}\) attain their norms form a dense subset of \(\mathcal P({}^n X)\). This resolves a problem posed by Aron, García and Maestre in 2003, extending their result for quadratic polynomials to arbitrary degrees without imposing any structural assumptions on the underlying Banach space. We further establish the corresponding theorem for vector-valued polynomials with arbitrary Banach space codomains. Our approach relies on a new perturbation scheme based on averaged linear contractions, which provides simultaneous control of the multilinear variables and yields norm attainment through iterated limits in the bidual.\end{abstract}
\maketitle
\enlargethispage{2pt}

\section{Introduction}

The study of norm attaining mappings begins with the Bishop--Phelps theorem \cite{BP}, which asserts that the norm attaining continuous linear functionals on a Banach space $X$ form a dense subset of $X^*$. The corresponding assertion for operators between arbitrary Banach spaces is false. Lindenstrauss \cite{L} showed, however, that every bounded linear operator can be approximated in norm by operators whose second adjoints attain their norms. Thus a density statement which fails on the original domain can be recovered by passing to the bidual.

The passage from linear functionals to homogeneous polynomials introduces a further difficulty. Norm attaining polynomials need not be dense in the space of continuous homogeneous polynomials; see Jim\'enez-Sevilla and Pay\'a \cite{JSP}. The situation changes when norm attainment is required of the canonical extension. Aron, Garc\'ia and Maestre \cite[Theorem 2.1]{AGM} proved that every continuous $2$-homogeneous scalar-valued polynomial can be approximated by polynomials whose canonical extensions to the bidual attain their norms. Their argument modifies the successive perturbations used by Lindenstrauss. In the same paper they pointed out that the analogue of Lindenstrauss' theorem for $n$-homogeneous polynomials with $n\ge3$ remained open \cite[p.~167]{AGM}.

For arbitrary degree, Carando, Lassalle and Mazzitelli \cite[Corollary 2.4]{CLM} obtained the scalar-valued polynomial Lindenstrauss theorem under the assumptions that $X^*$ is separable and has the approximation property. Their proof uses an integral representation associated with the symmetric tensor product, and their main result also covers polynomials with values in dual spaces. This settles the question for that class of spaces, but not for arbitrary Banach spaces. For vector-valued $2$-homogeneous polynomials, the analogue of \cite[Theorem 2.1]{AGM} was considered by Choi, Lee and Song \cite{CLS}. The goal of this note is to show that the scalar-valued polynomial Lindenstrauss theorem holds for every Banach space and in every degree (Theorem~\ref{thm:main}). In Section~\ref{sec:vector} we show that a modification of the argument also gives the vector-valued version, for arbitrary Banach spaces $X$ and $Y$ (Theorem~\ref{thm:vector}).

Throughout the paper, $\K$ denotes the scalar field, which is either $\R$ or $\C$; we write $\Rea a=a$ when $a\in\R$, and $X$ denotes a Banach space over $\K$. We write $B_X$ and $S_X$ for the closed unit ball and the unit sphere of $X$. For a positive integer $n$, a map $P:X\to\K$ is a continuous $n$-homogeneous polynomial if $P(x)=L(x,\ldots,x)$ for some continuous $n$-linear form $L:X\times\cdots\times X\to\K$, and $\Pn$ denotes the space of all such polynomials. Each $P\in\Pn$ has a unique associated continuous symmetric $n$-linear form $A=\check P$ (the symmetrization of $L$), with $P(x)=A(x,\ldots,x)$. We use the norms
\[
 \norm P=\sup_{x\in B_X}|P(x)|,
 \qquad
 \norm A=\sup_{x_1,\ldots,x_n\in B_X}|A(x_1,\ldots,x_n)|.
\]
By homogeneity and multilinearity these norms control all values:
\[
 |P(x)|\le\norm P\norm x^n,
 \qquad
 |A(x_1,\ldots,x_n)|\le\norm A\prod_{j=1}^n\norm{x_j}
 \qquad(x,x_1,\ldots,x_n\in X);
\]
for instance, if $x\ne0$ it suffices to apply the definition of $\norm P$ at $x/\norm x$. For $c_n=n^n/n!$, polarization gives
\begin{equation}\label{eq:polarnorm}
 \norm P\le\norm{\check P}\le c_n\norm P.
\end{equation}
For $x,u\in X$ and $0\le r\le n$ we write $A(x^{n-r},u^r)$ for the value of $A$ at the $n$-tuple formed by $n-r$ copies of $x$ and $r$ copies of $u$; by symmetry, their order is irrelevant.

For later use we recall the two polarization identities, valid over both fields: if $P\in\Pn$, $A=\check P$ and $v_1,\ldots,v_n\in X$, then
\begin{align}
 A(v_1,\ldots,v_n)
 &=\frac1{2^n n!}\sum_{\epsilon_1,\ldots,\epsilon_n=\pm1}
 \left(\prod_{j=1}^n\epsilon_j\right)
 P\left(\sum_{j=1}^n\epsilon_jv_j\right),\label{eq:signpolar}\\
 n!A(v_1,\ldots,v_n)
 &=\sum_{F\subseteq\{1,\ldots,n\}}(-1)^{n-|F|}
 P\left(\sum_{j\in F}v_j\right).\label{eq:finitepolar}
\end{align}
These identities, the estimate \eqref{eq:polarnorm} and the fact that $\Pn$ is a Banach space are standard; see, for instance, \cite{Di}. Moreover, replacing the arguments of $A(x,\ldots,x)$ by $y$ one at a time gives $P(x)-P(y)=\sum_{r=1}^nA(y^{r-1},x-y,x^{n-r})$, where the $r$-th term has $y$ in its first $r-1$ arguments, $x-y$ in the $r$-th one and $x$ in the remaining ones. If $x,y\in B_X$, each term has modulus at most $\norm{\check P}\norm{x-y}$, so \eqref{eq:polarnorm} yields the Lipschitz estimate
\begin{equation}\label{eq:lipschitz}
 |P(x)-P(y)|\le nc_n\norm P\,\norm{x-y}\qquad(x,y\in B_X).
\end{equation}

The multilinear density results established in \cite{AcGM} do not by themselves settle the polynomial problem: norm attainment of an extended form may occur at different points in its different variables, whereas the canonical extension of a polynomial must attain its norm at a single point $z\in B_{X^{**}}$, that is, at a point $(z,\ldots,z)$ of the diagonal. The proof of the quadratic theorem provided in \cite{AGM} controls mixed values at two independently chosen points. In higher degree, estimates with several arguments constrained to be equal do not control the successive limits defining the canonical extension.
For instance, when $n=3$ and $z\in X^{**}$ is the weak-star limit of a sequence $(x_k)\subseteq B_X$ along an ultrafilter, as in \eqref{eq:z} below, $\widetilde P(z)$ is an iterated limit of the values $\check P(x_i,x_j,x_k)$ in which the three indices vary independently (see \eqref{eq:ulpartial}). Lower bounds are therefore needed for all triples $i<j<k$, and not only for values such as $\check P(x_j,x_k,x_k)$.

To prove Theorem~\ref{thm:main}, instead of the perturbations used in \cite{AGM}, at each stage we perturb by an average of substitutions of the form
\[
 x\longmapsto (1-t)x+tf(x)u,
 \qquad \norm u=\norm f=f(u)=1,
 \qquad 0\le t\le1.
\]
These substitutions map the unit ball into itself and fix $u$. Perturbing by small multiples of their averages leads to an exact identity for certain nonnegative deficits; see \eqref{eq:deficitidentity}. Iterating this identity controls ordered compositions of substitutions, and a geometric lemma (Lemma~\ref{lem:geometry}) then gives lower bounds for the associated symmetric form at $n$ independently indexed arguments.

\section{The scalar-valued case}

Let us start by recalling the canonical extension. The dual $X^*$ consists of the continuous $\K$-linear functionals on $X$, and the bidual $X^{**}=(X^*)^*$ is also taken over $\K$. We identify $X$ with its canonical image under the linear isometry $J_X:X\to X^{**}$, $J_Xx(f)=f(x)$. The canonical extension of a continuous $n$-linear form $A$ on $X$ is obtained by extending first in the last variable and then successively in the preceding ones:
\begin{equation}\label{eq:arens}
 \overline A(z_1,\ldots,z_n)
 =z_1\bigl(x_1\mapsto z_2\bigl(x_2\mapsto\cdots
 z_n(x_n\mapsto A(x_1,\ldots,x_n))\cdots\bigr)\bigr).
\end{equation}
Here, for fixed $x_1,\ldots,x_{n-1}\in X$, the map $x_n\mapsto A(x_1,\ldots,x_n)$ belongs to $X^*$, and $z_n$ is applied to it. The resulting scalar is an $(n-1)$-linear function of $(x_1,\ldots,x_{n-1})$, bounded by $\norm A\norm{z_n}\prod_{j<n}\norm{x_j}$, and the same operation is repeated with $z_{n-1},\ldots,z_1$. Hence $\overline A$ is $n$-linear and $|\overline A(z_1,\ldots,z_n)|\le\norm A\prod_j\norm{z_j}$. Since $\overline A(J_Xx_1,\ldots,J_Xx_n)=A(x_1,\ldots,x_n)$, the form $\overline A$ extends $A$ and $\norm{\overline A}=\norm A$. If $A=\check P$, the canonical Aron--Berner extension \cite{AB} is
\begin{equation}\label{eq:AB}
 \widetilde P(z)=\overline A(z,\ldots,z),\qquad z\in X^{**}.
\end{equation}
Changing the order of extension of a symmetric form permutes the arguments of \eqref{eq:arens}, as one sees by relabelling the variables of $A$, so it does not change its diagonal. Although $\overline A$ need not be symmetric, its symmetrization has the same diagonal; hence \eqref{eq:AB} defines a continuous $n$-homogeneous polynomial on $X^{**}$.

We shall use \eqref{eq:arens} on the diagonal in the following recursive form. Fix $z\in X^{**}$ and a continuous symmetric $n$-linear form $A$ on $X$. Put $A_n=A$, and define
\begin{equation}\label{eq:partial}
 A_{r-1}(v_1,\ldots,v_{r-1})
 =z\bigl(v\mapsto A_r(v_1,\ldots,v_{r-1},v)\bigr),
 \quad r=n,\ldots,1.
\end{equation}
By downward induction each $A_r$ is a continuous $r$-linear form on $X$ with $\norm{A_r}\le\norm A\norm z^{n-r}$: for fixed $v_1,\ldots,v_{r-1}$, the functional $g(v)=A_r(v_1,\ldots,v_{r-1},v)$ belongs to $X^*$ and $\norm g\le\norm{A_r}\prod_{j<r}\norm{v_j}$, so $|A_{r-1}(v_1,\ldots,v_{r-1})|=|z(g)|\le\norm z\norm{A_r}\prod_{j<r}\norm{v_j}$; and $A_{r-1}$ is linear in each $v_j$, because $g$ depends linearly on each $v_j$ and $z$ is linear. Comparing with \eqref{eq:arens}, $A_0=\overline A(z,\ldots,z)$; in particular $A_0=\widetilde P(z)$ when $A=\check P$.

Our main result is the following.

\begin{samepage}
\begin{theorem}\label{thm:main}
Let $X$ be a real or complex Banach space and let $n\ge1$. The set of all $P\in\Pn$ whose canonical extensions to $X^{**}$ attain their norms is dense in $\Pn$. More precisely, for every $P\in\Pn$ and every $\eps>0$ there exist $Q\in\Pn$ and $z\in B_{X^{**}}$ such that
\[
 \norm{Q-P}<\eps,
 \qquad
 |\widetilde Q(z)|=\norm{\widetilde Q}=\norm Q.
\]
\end{theorem}
\end{samepage}

The proof is organized as follows. Lemma~\ref{lem:isometry} recalls that the canonical extension preserves norms. Lemma~\ref{lem:K} is an elementary estimate for nonnegative real polynomials of one real variable, which turns a bound for an integral into a bound for every value of the integrand. After introducing the substitutions $S_{u,f}(t)$ and their averages $\Rop_{u,f}$, Lemma~\ref{lem:geometry} shows that a polynomial whose real part is almost maximal along all compositions of $n-1$ such substitutions applied to a fixed vector has an almost maximal fully mixed value of its symmetric form. The polynomial $Q$ of Theorem~\ref{thm:main} is the limit of successive perturbations $P_{k+1}=P_k+\alpha_k\Rop_kP_k$; see \eqref{eq:iteration}. An exact identity for the resulting deficits yields the hypothesis of Lemma~\ref{lem:geometry} for every increasing $n$-tuple of almost norming points, and iterated limits along an ultrafilter turn the resulting lower bounds into $|\widetilde Q(z)|\ge\norm Q$. The vector-valued case is treated in Section~\ref{sec:vector}.

We give the preliminary extension facts and estimates before proving the theorem. The first lemma is the norm preservation theorem for canonical polynomial extensions, due to Davie and Gamelin \cite{DG}.

\begin{lemma}[{\cite[Theorem~3]{DG}}]\label{lem:isometry}
For every $P\in\Pn$, the canonical extension satisfies $\norm{\widetilde P}=\norm P$. The map $P\mapsto\widetilde P$ is linear.
\end{lemma}

Davie and Gamelin derive \cite[Theorem~3]{DG} from \cite[Theorem~1]{DG}, the approximation of points of $B_{X^{**}}$ by arithmetic means of points of $B_X$ in the polynomial-star topology; the canonical extension in \cite{DG} is the one in \eqref{eq:arens}, obtained by extending the last variable first. Although \cite{DG} is written for complex Banach spaces, the proof of \cite[Theorem~1]{DG} (the lemma in \cite[\S3]{DG} and the averaging argument that follows it) uses only the multilinearity of the forms, the weak-star continuity of the extended forms in each variable when the preceding variables lie in $X$, and Goldstine's theorem; therefore it also applies to real Banach spaces. Only the inequality $\norm{\widetilde P}\le\norm P$ is nontrivial: the reverse inequality $\norm{\widetilde P}\ge\norm P$ holds because $\widetilde P(J_Xx)=P(x)$ and $J_X$ maps $B_X$ into $B_{X^{**}}$. Linearity of $P\mapsto\widetilde P$ follows from \eqref{eq:arens} and \eqref{eq:AB}, since $P\mapsto\check P$ is linear and $\overline A$ depends linearly on $A$.

The next estimate is a finite-dimensional fact about ordinary real polynomials.

\begin{lemma}\label{lem:K}
For every $n\ge1$ there is a constant $K_n\ge1$ such that
\begin{equation}\label{eq:K}
 \max_{0\le t\le1}r(t)\le K_n\int_0^1r(t)\,dt
\end{equation}
for every real polynomial $r$ of degree at most $n$ which is nonnegative on $[0,1]$.
\end{lemma}

\begin{proof}
Assume that this is not the case. Then, for every $j\in\N$, there is a real polynomial $r_j$ of degree at most $n$, nonnegative on $[0,1]$, such that $\max_{[0,1]}r_j>j\int_0^1r_j(t)\,dt$. Such an $r_j$ is not identically zero, so $\norm{r_j}_\infty=\max_{[0,1]}r_j>0$, where $\norm\cdot_\infty$ denotes the maximum norm on $[0,1]$. Replacing $r_j$ by $r_j/\norm{r_j}_\infty$, we may assume that
\[
 \norm{r_j}_\infty=1
 \qquad\text{and}\qquad
 \int_0^1r_j(t)\,dt<\frac1j .
\]
Fix distinct points $t_0,\ldots,t_n\in[0,1]$ and write $r_j(t)=\sum_{k=0}^na_{j,k}t^k$. The coefficient vector $a_j=(a_{j,0},\ldots,a_{j,n})$ solves the linear system
\[
 \sum_{k=0}^na_{j,k}t_i^k=r_j(t_i)\qquad(0\le i\le n),
\]
whose matrix $(t_i^k)_{i,k}$ is invertible, since a polynomial of degree at most $n$ with $n+1$ distinct zeros vanishes identically. As $|r_j(t_i)|\le1$, the vectors $a_j$ form a bounded sequence in $\R^{n+1}$, and a subsequence $(a_{j_m})_m$ converges to some $a=(a_0,\ldots,a_n)$. Put $r(t)=\sum_{k=0}^na_kt^k$. For $t\in[0,1]$ we have $|r_{j_m}(t)-r(t)|\le\sum_{k=0}^n|a_{j_m,k}-a_k|$, so $r_{j_m}\to r$ uniformly on $[0,1]$. Consequently $r\ge0$ on $[0,1]$, $\norm r_\infty=\lim_m\norm{r_{j_m}}_\infty=1$, and
\[
 \int_0^1r(t)\,dt=\lim_{m\to\infty}\int_0^1r_{j_m}(t)\,dt=0 .
\]
A continuous nonnegative function with zero integral vanishes on $[0,1]$: if $r(s)>0$ for some $s\in[0,1]$, then $r>r(s)/2$ on an interval of positive length containing $s$, and the integral of $r$ would be positive. Thus $r=0$, which contradicts $\norm r_\infty=1$.
\end{proof}

The constant $K_n$ depends only on $n$. This will allow us to use the same constant at every stage of the construction, since the perturbations never increase the degree.

For $u\in S_X$ and $f\in S_{X^*}$ with $f(u)=1$, and for $t\in[0,1]$, define the map $S_{u,f}(t):X\to X$ and the averaging operator $\Rop_{u,f}$ on $\Pn$ by
\begin{equation}\label{eq:substitution}
 S_{u,f}(t)x=(1-t)x+tf(x)u\quad(x\in X),
 \qquad
 \Rop_{u,f}P(x)=\int_0^1P(S_{u,f}(t)x)\,dt\quad(x\in X).
\end{equation}
The integral is that of a continuous scalar function of the real variable $t$; in fact $t\mapsto P(S_{u,f}(t)x)$ is a polynomial in $t$, see \eqref{eq:expansion} below.

Each $S_{u,f}(t)$ is a $\K$-linear operator on $X$, because $f$ is $\K$-linear; the parameter $t$ is real even when $\K=\C$. It is a contraction, that is, $\norm{S_{u,f}(t)x}\le\norm x$ for every $x\in X$: since $\norm u=\norm f=1$ and $t,1-t\ge0$,
\[
 \norm{S_{u,f}(t)x}\le(1-t)\norm x+t|f(x)|\norm u\le(1-t)\norm x+t\norm x=\norm x .
\]
In particular $S_{u,f}(t)$ maps $B_X$ into itself. Moreover $S_{u,f}(t)u=(1-t)u+tf(u)u=u$, so every $S_{u,f}(t)$ fixes $u$. We show below that $\Rop_{u,f}P\in\Pn$. Clearly $\Rop_{u,f}$ is linear, and integrating over $t\in[0,1]$ the inequality $|P(S_{u,f}(t)x)|\le\norm P$ $(x\in B_X)$ and the identity $P(S_{u,f}(t)u)=P(u)$, we obtain
\begin{equation}\label{eq:Rproperties}
 \norm{\Rop_{u,f}P}\le\norm P,
 \qquad
 \Rop_{u,f}P(u)=P(u).
\end{equation}

To see that $\Rop_{u,f}P$ belongs to $\Pn$, write $A=\check P$. Expand each of the $n$ arguments of $A((1-t)x+tf(x)u,\ldots,(1-t)x+tf(x)u)$ as a sum of two terms. The $\binom nr$ resulting terms in which exactly $r$ arguments equal $tf(x)u$ coincide by symmetry, so
\begin{equation}\label{eq:expansion}
 P(S_{u,f}(t)x)=\sum_{r=0}^n\binom nr(1-t)^{n-r}t^rf(x)^rA(x^{n-r},u^r).
\end{equation}
For nonnegative integers $r,s$ we have $\int_0^1(1-t)^s\,dt=1/(s+1)$ and, if $r\ge1$, integration by parts gives $\int_0^1t^r(1-t)^s\,dt=\frac r{s+1}\int_0^1t^{r-1}(1-t)^{s+1}\,dt$. By induction on $r$,
\[
 \int_0^1t^r(1-t)^s\,dt=\frac{r!\,s!}{(r+s+1)!}.
\]
With $s=n-r$ this gives
\[
 \binom nr\int_0^1t^r(1-t)^{n-r}\,dt=\frac{n!}{r!\,(n-r)!}\cdot\frac{r!\,(n-r)!}{(n+1)!}=\frac1{n+1},
\]
and integrating \eqref{eq:expansion} over $[0,1]$ yields
\begin{equation}\label{eq:Rfinite}
 \Rop_{u,f}P(x)=\frac1{n+1}\sum_{r=0}^n
 f(x)^r\check P(x^{n-r},u^r).
\end{equation}
Each summand is a continuous $n$-homogeneous polynomial in $x$: it is the restriction to the diagonal of the continuous $n$-linear form
\[
 (x_1,\ldots,x_n)\longmapsto f(x_1)\cdots f(x_r)\,\check P(x_{r+1},\ldots,x_n,u,\ldots,u)
\]
with $r$ copies of $u$, and the symmetrization of this form has the same diagonal. Hence $\Rop_{u,f}P\in\Pn$.

Before we continue, let us point out that averages of different substitutions compose in reverse order. If $\Rop_1=\Rop_{u_1,f_1}$ and $\Rop_2=\Rop_{u_2,f_2}$, then, by \eqref{eq:substitution} applied twice,
\[
 (\Rop_2\Rop_1P)(x)=\int_0^1(\Rop_1P)(S_{u_2,f_2}(t)x)\,dt
 =\int_0^1\!\!\int_0^1P\bigl(S_{u_1,f_1}(s)S_{u_2,f_2}(t)x\bigr)\,ds\,dt .
\]
Thus the substitution of the average applied first to $P$ acts last on the vector. Ordered compositions of this kind appear in the proof of Theorem~\ref{thm:main}.

We now give the estimate which allows the substitutions to control all the variables of the associated symmetric form. For a fixed $n\ge2$, let
\begin{equation}\label{eq:constants}
 D_n=\frac1{n!}\sum_{r=1}^n\binom nr r^n,
 \qquad
 H_n=D_n\bigl(\sqrt2+2n(n-1)c_n\bigr),
\end{equation}
where, as above, $c_n=n^n/n!$. Both constants depend only on $n$, and $H_n\ge D_n\ge1$, because the term $r=n$ gives $D_n\ge n^n/n!\ge1$.

\begin{lemma}\label{lem:geometry}
Let $P\in\Pn$, let $M>0$ be such that $\norm P\le M$, and take $u_1,\ldots,u_n\in B_X$ and $f_1,\ldots,f_{n-1}\in B_{X^*}$. For $1\le i<n$ and $t\in[0,1]$, let $S_i(t):X\to X$ be the linear map
\[
 S_i(t)x=(1-t)x+tf_i(x)u_i\qquad(x\in X).
\]
Suppose that, for some $E\ge0$,
\begin{equation}\label{eq:path}
 \Rea P\bigl(S_1(t_1)\cdots S_{n-1}(t_{n-1})u_n\bigr)\ge M-E
 \quad\text{for every } t_1,\ldots,t_{n-1}\in[0,1].
\end{equation}
Then
\begin{equation}\label{eq:geometry}
 |\check P(u_1,\ldots,u_n)|\ge M-H_n\sqrt{ME}.
\end{equation}
\end{lemma}

In \eqref{eq:path} the operators act from right to left: $S_{n-1}(t_{n-1})$ is applied first to $u_n$, and $S_1(t_1)$ last. Since $\norm{f_i}\le1$ and $\norm{u_i}\le1$, the estimate used above for $S_{u,f}(t)$ shows that every $S_i(t)$ is a linear contraction.

\begin{proof}
If $E\ge M$, then $\sqrt{ME}\ge M$ and $H_n\ge1$, so the right-hand side of \eqref{eq:geometry} is at most $M-M=0$ and the assertion is immediate. Assume $0\le E<M$ and put $\rho=E/M\in[0,1)$. Set
\[
 b_i=f_i(u_n)\quad(i<n),\qquad b_n=1,\qquad z_i=b_i u_i\quad(i\le n).
\]
Since $|b_i|\le\norm{f_i}\norm{u_n}\le1$, every $z_i$ belongs to $B_X$; note that $z_n=u_n$.

To start, we show that the functionals $f_i$ are almost constant on the later points $z_j$, that is, $f_i(z_j)$ is close to $b_i$ whenever $i<j\le n$. For this, fix $i<j<n$. Choose $t_i=1$, $t_j=1/2$, and all the other parameters zero in \eqref{eq:path}. We obtain
\[
 \Rea\left[\left(\frac{b_i+b_jf_i(u_j)}2\right)^nP(u_i)\right]\ge M-E.
\]
Indeed, $S_h(0)$ is the identity for every $h$, so the composition in \eqref{eq:path} reduces to $S_i(1)S_j(1/2)u_n$. First,
\[
 S_j(1/2)u_n=\tfrac12u_n+\tfrac12f_j(u_n)u_j=\tfrac12u_n+\tfrac12b_ju_j .
\]
Next, $S_i(1)v=f_i(v)u_i$ for every $v\in X$; hence, by the linearity of $f_i$,
\[
 S_i(1)S_j(1/2)u_n=f_i\bigl(\tfrac12u_n+\tfrac12b_ju_j\bigr)u_i
 =\tfrac12\bigl(b_i+b_jf_i(u_j)\bigr)u_i ,
\]
and the displayed inequality follows from the homogeneity $P(cu_i)=c^nP(u_i)$, $c\in\K$. Since the real part of a scalar is at most its modulus, and $|P(u_i)|\le M$ because $\norm P\le M$ and $u_i\in B_X$, this inequality gives $M\bigl|(b_i+b_jf_i(u_j))/2\bigr|^n\ge M-E$. Dividing by $M>0$, we obtain
\[
 \left|\frac{b_i+b_jf_i(u_j)}2\right|^n\ge1-\rho .
\]

Both scalars $\xi=b_i$ and $\upsilon=b_jf_i(u_j)$ have modulus at most one, and the previous inequality says that $|\xi+\upsilon|^2\ge4(1-\rho)^{2/n}$. By the parallelogram identity in $\K$,
\[
 \begin{split}
 |\xi-\upsilon|^2
 &=2(|\xi|^2+|\upsilon|^2)-|\xi+\upsilon|^2\\
 &\le4-4(1-\rho)^{2/n}\le4\rho,
 \end{split}
\]
where the last inequality holds because $0<1-\rho\le1$ and $2/n\le1$ (recall that $n\ge2$), so that $(1-\rho)^{2/n}\ge1-\rho$. Since $f_i(z_j)=b_jf_i(u_j)=\upsilon$ and $b_i=\xi$, and since $f_i(z_n)=f_i(u_n)=b_i$, it follows that
\begin{equation}\label{eq:compatibility}
 |f_i(z_j)-b_i|\le2\sqrt\rho\quad(i<j<n),
 \qquad f_i(z_n)=b_i\quad(i<n).
\end{equation}
When $n=2$ there are no pairs $i<j<n$, and \eqref{eq:compatibility} reduces to the identity $f_1(z_2)=b_1$.

Next, we compare the composed vectors appearing in \eqref{eq:path} with convex combinations of $z_1,\ldots,z_n$. Fix $t_1,\ldots,t_{n-1}\in[0,1]$ and define recursively, starting from the right,
\[
 y_n=v_n=u_n,\qquad
 y_i=S_i(t_i)y_{i+1},\qquad
 v_i=(1-t_i)v_{i+1}+t_i z_i\qquad(i=n-1,\ldots,1).
\]
Thus $y_1=S_1(t_1)\cdots S_{n-1}(t_{n-1})u_n$ is the vector appearing in \eqref{eq:path}, and $y_1\in B_X$ because each $S_i(t_i)$ is a contraction and $u_n\in B_X$. Since $v_n=z_n$ and each $v_i$ is a convex combination of $z_i$ and $v_{i+1}$, downward induction shows that $v_{i+1}$ is a convex combination of $z_{i+1},\ldots,z_n$ for every $i<n$; in particular $v_1\in\conv\{z_1,\ldots,z_n\}\subseteq B_X$. Subtracting the definitions of $y_i$ and $v_i$, using $z_i=b_iu_i$, and writing $f_i(y_{i+1})=f_i(y_{i+1}-v_{i+1})+f_i(v_{i+1})$, we get the identity
\[
 \begin{split}
 y_i-v_i&=S_i(t_i)y_{i+1}-(1-t_i)v_{i+1}-t_iz_i\\
 &=(1-t_i)y_{i+1}+t_if_i(y_{i+1})u_i-(1-t_i)v_{i+1}-t_ib_iu_i\\
 &=(1-t_i)(y_{i+1}-v_{i+1})+t_i\bigl(f_i(y_{i+1})-b_i\bigr)u_i\\
 &=(1-t_i)(y_{i+1}-v_{i+1})+t_if_i(y_{i+1}-v_{i+1})u_i+t_i\bigl(f_i(v_{i+1})-b_i\bigr)u_i\\
 &=S_i(t_i)(y_{i+1}-v_{i+1})+t_i\bigl(f_i(v_{i+1})-b_i\bigr)u_i .
 \end{split}
\]
Writing $v_{i+1}=\sum_{j>i}\lambda_jz_j$ with $\lambda_j\ge0$ and $\sum_{j>i}\lambda_j=1$, we get
\[
 f_i(v_{i+1})-b_i=\sum_{j>i}\lambda_j\bigl(f_i(z_j)-b_i\bigr),
\]
which has modulus at most $2\sqrt\rho$ by \eqref{eq:compatibility}. Since $S_i(t_i)$ is a contraction, $0\le t_i\le1$ and $\norm{u_i}\le1$, we obtain
\[
 \norm{y_i-v_i}\le\norm{y_{i+1}-v_{i+1}}+2\sqrt\rho .
\]
As $y_n-v_n=0$, the $n-1$ steps $i=n-1,\ldots,1$ give
\begin{equation}\label{eq:convexdistance}
 \norm{y_1-v_1}\le2(n-1)\sqrt\rho.
\end{equation}

As the parameters vary, $v_1$ runs through the whole convex hull $\conv\{z_1,\ldots,z_n\}$. To see this, let $\theta_1,\ldots,\theta_n\ge0$ with $\sum_i\theta_i=1$; we choose parameters $t_1,\ldots,t_{n-1}\in[0,1]$ for which $v_1=\sum_{i=1}^n\theta_iz_i$. For $1\le i\le n$ let $r_i=\sum_{h\ge i}\theta_h=1-\sum_{h<i}\theta_h$ be the total weight of $\theta_i,\ldots,\theta_n$, so that $r_1=1$, $r_n=\theta_n$, $0\le\theta_i\le r_i$, and $r_{i+1}=r_i-\theta_i$ for $i<n$. For $i=1,\ldots,n-1$ put
\[
 t_i=\theta_i/r_i\ \text{ if } r_i>0,
 \qquad
 t_i=0\ \text{ if } r_i=0 .
\]
Then $t_i\in[0,1]$, and $r_{i+1}=r_i(1-t_i)$ in both cases (if $r_i=0$, then $\theta_i=\cdots=\theta_n=0$ and $r_{i+1}=0$). Hence $r_i=\prod_{h<i}(1-t_h)$ for every $i\le n$. Since $\theta_i=t_ir_i$ for $i<n$ (if $r_i=0$, both sides vanish), we obtain
\begin{equation}\label{eq:weights}
 \theta_i=t_ir_i=t_i\prod_{h<i}(1-t_h)\quad(i<n),
 \qquad
 \theta_n=r_n=\prod_{h<n}(1-t_h).
\end{equation}
On the other hand, for any choice of the parameters, downward induction on $i$ shows that
\[
 v_i=\sum_{h=i}^{n-1}t_h\prod_{g=i}^{h-1}(1-t_g)\,z_h+\prod_{g=i}^{n-1}(1-t_g)\,z_n
 \qquad(1\le i\le n),
\]
since $v_n=z_n$ and $v_i=t_iz_i+(1-t_i)v_{i+1}$; for $i=1$ this is $v_1=\sum_{i=1}^n\theta_iz_i$ with precisely the weights in \eqref{eq:weights}. Thus every point of the convex hull, including the points of its faces such as the averages of subsets of $\{z_1,\ldots,z_n\}$, is of the form $v_1$ for suitable parameters.

Let us now show that $P$ is almost constant, with value close to $M$, on the convex hull. For every choice of the parameters, \eqref{eq:path} gives $\Rea P(y_1)\ge M-E$, and $|P(y_1)|\le M$ because $\norm P\le M$ and $y_1\in B_X$. Hence
\[
 |P(y_1)-M|^2
 =|P(y_1)|^2+M^2-2M\Rea P(y_1)\le M^2+M^2-2M(M-E)=2ME.
\]
Now let $v\in\conv\{z_1,\ldots,z_n\}$. By the previous paragraph, there are parameters $t_1,\ldots,t_{n-1}\in[0,1]$ for which the vector $v_1$ defined by the recursion above equals $v$; for these same parameters, let $y_1=S_1(t_1)\cdots S_{n-1}(t_{n-1})u_n$. Both $v$ and $y_1$ lie in $B_X$, so \eqref{eq:lipschitz} (with $\norm P\le M$), \eqref{eq:convexdistance} and the previous inequality give
\[
 \begin{split}
 |P(v)-M|
 &\le|P(v)-P(y_1)|+|P(y_1)-M|\\
 &\le nc_nM\cdot2(n-1)\sqrt{E/M}+\sqrt{2ME}\\
 &=\bigl(2n(n-1)c_n+\sqrt2\bigr)\sqrt{ME}.
 \end{split}
\]
That is,
\begin{equation}\label{eq:flat}
 |P(v)-M|\le C'_n\sqrt{ME}
 \quad\bigl(v\in\conv\{z_1,\ldots,z_n\}\bigr),
 \qquad C'_n=\sqrt2+2n(n-1)c_n.
\end{equation}

Finally, we estimate the fully mixed value $\check P(z_1,\ldots,z_n)$. Let $F\subseteq\{1,\ldots,n\}$ be nonempty, with $|F|=r$. The average $r^{-1}\sum_{j\in F}z_j$ belongs to $\conv\{z_1,\ldots,z_n\}$, and $P(\sum_{j\in F}z_j)=r^nP(r^{-1}\sum_{j\in F}z_j)$ by homogeneity. Therefore \eqref{eq:flat} gives
\[
 \left|P\left(\sum_{j\in F}z_j\right)-r^nM\right|
 =r^n\left|P\left(\frac1r\sum_{j\in F}z_j\right)-M\right|
 \le r^nC'_n\sqrt{ME}.
\]
Applying \eqref{eq:finitepolar} to the polynomial $t\mapsto t^n$ on $\K$, whose associated form is $(t_1,\ldots,t_n)\mapsto t_1\cdots t_n$, with $v_1=\cdots=v_n=1$, we obtain
\[
 \sum_{F\subseteq\{1,\ldots,n\}}(-1)^{n-|F|}|F|^n=n!.
\]
Subtracting $M$ times this identity from \eqref{eq:finitepolar} applied to $z_1,\ldots,z_n$, we get
\[
 n!\bigl(\check P(z_1,\ldots,z_n)-M\bigr)
 =\sum_{F\subseteq\{1,\ldots,n\}}(-1)^{n-|F|}
 \left(P\left(\sum_{j\in F}z_j\right)-|F|^nM\right).
\]
The term corresponding to $F=\emptyset$ vanishes, because $P(0)=0$ and $|\emptyset|^nM=0$. Since there are $\binom nr$ subsets of size $r$, the triangle inequality and the previous estimate give
\[
 \begin{split}
 n!\,|\check P(z_1,\ldots,z_n)-M|
 &\le\sum_{\emptyset\ne F\subseteq\{1,\ldots,n\}}\left|P\left(\sum_{j\in F}z_j\right)-|F|^nM\right|\\
 &\le\sum_{r=1}^n\binom nr r^nC'_n\sqrt{ME}=n!\,D_nC'_n\sqrt{ME},
 \end{split}
\]
where the last equality is the definition of $D_n$ in \eqref{eq:constants}; that is, $|\check P(z_1,\ldots,z_n)-M|\le D_nC'_n\sqrt{ME}=H_n\sqrt{ME}$. It remains to compare with $\check P(u_1,\ldots,u_n)$: by multilinearity and $b_n=1$,
\[
 \check P(z_1,\ldots,z_n)
 =\left(\prod_{i<n}b_i\right)\check P(u_1,\ldots,u_n),
 \qquad |b_i|\le1.
\]
Therefore
\[
 |\check P(u_1,\ldots,u_n)|\ge|\check P(z_1,\ldots,z_n)|\ge M-H_n\sqrt{ME},
\]
which proves the lemma.
\end{proof}

We now prove Theorem~\ref{thm:main}. The proof is based on a careful perturbation of $P$ by averages of the substitutions $S_{u,f}(t)$.

\begin{proof}[Proof of Theorem~\ref{thm:main}]
Let $P\in\Pn$ and $\eps>0$. The case $n=1$ is trivial: if $P=f\in X^*$ and $f\ne0$, then $\widetilde f(z)=z(f)$ by \eqref{eq:arens}, and the Hahn--Banach theorem, applied in $X^*$, provides $z\in B_{X^{**}}$ with $z(f)=\norm f=\norm{\widetilde f}$, so we may take $Q=P$. If $P=0$, we take $Q=0$. Moreover, if the theorem holds for $P/\norm P$ and $\eps/\norm P$, with a polynomial $R$ and a point $z$, then it holds for $P$ and $\eps$ with $Q=\norm P\,R$ and the same $z$, by the linearity of $P\mapsto\widetilde P$. Hence, without loss of generality, we may assume that $n\ge2$ and $\norm P=1$. Choose
\begin{equation}\label{eq:alphas}
 0<a<\min\{1/8,\eps/4\},\qquad\text{and define}\qquad
 \alpha_k=a^{2^{k-1}}\quad(k\ge1).
\end{equation}
Then $\alpha_{k+1}=\alpha_k^2$ and $0<\alpha_k\le a<1/8$ for every $k$, and
\begin{equation}\label{eq:tails}
 \sum_{h=k}^{\infty}\alpha_h\le2\alpha_k,
 \qquad
 \sum_{h=k}^{\infty}\alpha_h^2\le2\alpha_k^2.
\end{equation}
Indeed, for any nonnegative integer $r$ we have $2^r\ge r+1$, hence $\alpha_{k+r}=\alpha_k^{2^r}\le\alpha_k^{r+1}$ and $\sum_{h\ge k}\alpha_h\le\sum_{r\ge0}\alpha_k^{r+1}=\alpha_k/(1-\alpha_k)\le2\alpha_k$; moreover $\sum_{h\ge k}\alpha_h^2=\sum_{h\ge k+1}\alpha_h\le2\alpha_{k+1}=2\alpha_k^2$. The fast decay $\alpha_{k+1}=\alpha_k^2$ will be essential in the proof of \eqref{eq:compositions} below.

Starting with $P_1=P$, we construct inductively polynomials $P_k\in\Pn$, vectors $x_k\in S_X$ and functionals $f_k\in S_{X^*}$ such that, with $M_k=\norm{P_k}$,
\begin{equation}\label{eq:nearmax}
 |P_k(x_k)|\ge M_k-\alpha_k^2,
 \qquad f_k(x_k)=1,
\end{equation}
and
\begin{equation}\label{eq:iteration}
 P_{k+1}=P_k+\alpha_k\Rop_{x_k,f_k}P_k,
 \qquad A_k=\check P_k.
\end{equation}
For simplicity, we abbreviate $S_k(t)=S_{x_k,f_k}(t)$ and $\Rop_k=\Rop_{x_k,f_k}$. Once $P_k$ is known, the choice in \eqref{eq:nearmax} is possible: the supremum of $|P_k|$ over $B_X$ equals its supremum over $S_X$, because $|P_k(x/\norm x)|=\norm x^{-n}|P_k(x)|\ge|P_k(x)|$ when $0<\norm x\le1$, so some $x_k\in S_X$ satisfies $|P_k(x_k)|\ge M_k-\alpha_k^2$; and the Hahn--Banach theorem provides $f_k\in S_{X^*}$ with $f_k(x_k)=1$. By \eqref{eq:Rfinite}, $P_{k+1}\in\Pn$.

By \eqref{eq:Rproperties} and \eqref{eq:nearmax},
\begin{align}
 \norm{P_{k+1}-P_k}&\le\alpha_k M_k,\label{eq:increment}\\
 (1+\alpha_k)(M_k-\alpha_k^2)
 &\le M_{k+1}\le(1+\alpha_k)M_k.\label{eq:normbounds}
\end{align}
Indeed, by \eqref{eq:Rproperties}, $\norm{P_{k+1}-P_k}=\alpha_k\norm{\Rop_kP_k}\le\alpha_kM_k$, which is \eqref{eq:increment}; together with the triangle inequality, this gives the upper bound in \eqref{eq:normbounds}. For the lower bound, every $S_k(t)$ fixes $x_k$, so $\Rop_kP_k(x_k)=P_k(x_k)$ and $P_{k+1}(x_k)=(1+\alpha_k)P_k(x_k)$. Since $1+\alpha_k>0$ and $x_k\in B_X$, we get $M_{k+1}\ge|P_{k+1}(x_k)|=(1+\alpha_k)|P_k(x_k)|\ge(1+\alpha_k)(M_k-\alpha_k^2)$. If $M_k\ge1$, then
\[
 M_{k+1}-M_k\ge\alpha_kM_k-(1+\alpha_k)\alpha_k^2
 \ge\alpha_k(1-\alpha_k-\alpha_k^2)>0,
\]
because $\alpha_k<1/8$. Since $M_1=1$, it follows by induction that $(M_k)$ is increasing and $M_k\ge1$ for every $k$. On the other hand, iterating the upper bound in \eqref{eq:normbounds} from $M_1=1$ gives $M_k\le\prod_{h=1}^{k-1}(1+\alpha_h)$. Since $\alpha_h=a^{2^{h-1}}$, applying the identity $(1-b)(1+b)=1-b^2$ successively with $b=a,a^2,a^4,\ldots,a^{2^{k-2}}$ yields
\[
 (1-a)\prod_{h=1}^{k-1}\bigl(1+a^{2^{h-1}}\bigr)
 =(1-a^2)\prod_{h=2}^{k-1}\bigl(1+a^{2^{h-1}}\bigr)
 =\cdots
 =1-a^{2^{k-1}}
 =1-\alpha_k .
\]
Dividing by $1-a$, we obtain
\begin{equation}\label{eq:product}
 M_k\le\prod_{h<k}(1+\alpha_h)
 =\frac{1-\alpha_k}{1-a}<\frac1{1-a}<2.
\end{equation}
For $m>k$, \eqref{eq:increment}, \eqref{eq:product} and \eqref{eq:tails} give
\[
 \norm{P_m-P_k}\le\sum_{h=k}^{m-1}\alpha_hM_h\le\frac1{1-a}\sum_{h=k}^\infty\alpha_h\le\frac{2\alpha_k}{1-a}.
\]
Hence $(P_k)$ is a Cauchy sequence in the Banach space $\Pn$. Let $Q\in\Pn$ be its limit and $M=\norm Q$. Letting $m\to\infty$, and using $2/(1-a)<4$, we obtain
\begin{equation}\label{eq:tailQ}
 \norm{Q-P_k}\le4\alpha_k,
 \qquad
 \norm{Q-P}\le\frac{2a}{1-a}<4a<\eps.
\end{equation}
Since $|\norm Q-M_k|\le\norm{Q-P_k}$, we have $M_k\to M$; since $M_k\ge1$ for every $k$ and $M_k<(1-a)^{-1}$ by \eqref{eq:product}, we get $1\le M\le(1-a)^{-1}<2$. Moreover, $\check Q-A_k$ is the symmetric form associated with $Q-P_k$, so \eqref{eq:polarnorm} and \eqref{eq:tailQ} give
\[
 \norm{\check Q-A_k}\le c_n\norm{Q-P_k}\le4c_n\alpha_k .
\]

For $\omega\in\T_\K=\{\lambda\in\K:|\lambda|=1\}$ and $y\in B_X$, put
\begin{equation}\label{eq:deficits}
 d_k^\omega(y)=M_k-\Rea(\omega P_k(y)),
 \qquad
 \lambda_k=(1+\alpha_k)M_k-M_{k+1}.
\end{equation}
Over $\R$ we have $\T_\R=\{-1,1\}$. Then $d_k^\omega\ge0$ on $B_X$, since $\Rea(\omega P_k(y))\le|\omega P_k(y)|=|P_k(y)|\le M_k$. By \eqref{eq:normbounds}, $\lambda_k\ge0$ and $\lambda_k\le(1+\alpha_k)M_k-(1+\alpha_k)(M_k-\alpha_k^2)=(1+\alpha_k)\alpha_k^2$; together with \eqref{eq:tails}, this gives
\begin{equation}\label{eq:loss}
 0\le\lambda_k\le2\alpha_k^2,
 \qquad
 \sum_{h=k}^{\infty}\lambda_h\le4\alpha_k^2.
\end{equation}
Indeed, since $\alpha_h<1$,
\[
 \sum_{h=k}^\infty\lambda_h\le\sum_{h=k}^\infty(1+\alpha_h)\alpha_h^2\le\sum_{h=k}^\infty2\alpha_h^2\le4\alpha_k^2 ,
\]
where the last inequality is \eqref{eq:tails}.
For every fixed $\omega\in\T_\K$ and every $y\in B_X$, the recurrence \eqref{eq:iteration} gives the exact identity
\begin{equation}\label{eq:deficitidentity}
 \begin{split}
 &d_k^\omega(y)+\alpha_k\int_0^1d_k^\omega(S_k(t)y)\,dt\\
 &\qquad=(1+\alpha_k)M_k-\Rea\left(\omega\left[P_k(y)+\alpha_k\int_0^1P_k(S_k(t)y)\,dt\right]\right)\\
 &\qquad=(1+\alpha_k)M_k-\Rea\bigl(\omega P_{k+1}(y)\bigr)\\
 &\qquad=\bigl[M_{k+1}-\Rea(\omega P_{k+1}(y))\bigr]+\bigl[(1+\alpha_k)M_k-M_{k+1}\bigr]\\
 &\qquad=d_{k+1}^\omega(y)+\lambda_k.
 \end{split}
\end{equation}
In the first equality we used that $\int_0^1M_k\,dt=M_k$ and that the real part commutes with the integral; the second one is \eqref{eq:iteration}, and the last one is \eqref{eq:deficits}.

Since $y\in B_X$, we have $S_k(t)y\in B_X$ for every $t$, so both terms on the left-hand side of the previous equalities are nonnegative. Dropping the integral term, or dropping the first term, we obtain
\[
 d_k^\omega(y)\le d_{k+1}^\omega(y)+\lambda_k,
 \qquad
 \int_0^1d_k^\omega(S_k(t)y)\,dt\le\frac{d_{k+1}^\omega(y)+\lambda_k}{\alpha_k}.
\]
The function $r(t)=d_k^\omega(S_k(t)y)$ is a real polynomial of degree at most $n$ in the real variable $t$: by \eqref{eq:expansion}, $P_k(S_k(t)y)$ is a polynomial of degree at most $n$ in $t$ with coefficients in $\K$, and taking the real part of $\omega$ times it produces a real polynomial. Since $r\ge0$ on $[0,1]$, Lemma~\ref{lem:K} gives $\max_{[0,1]}r\le K_n\int_0^1r$. We have proved, for every $\omega\in\T_\K$, every $y\in B_X$ and every $k\ge1$,
\begin{align}
 d_k^\omega(y)&\le d_{k+1}^\omega(y)+\lambda_k,\label{eq:backward}\\
 \sup_{t\in[0,1]}d_k^\omega(S_k(t)y)
 &\le\frac{K_n}{\alpha_k}\bigl(d_{k+1}^\omega(y)+\lambda_k\bigr).
 \label{eq:insert}
\end{align}
The first inequality passes from stage $k+1$ to stage $k$ at the same point $y$. The second replaces $y$ by any of the points $S_k(t)y$, $t\in[0,1]$, at the cost of the factor $K_n/\alpha_k$.

For each $\ell\in\N$, choose
\begin{equation}\label{eq:orientation}
 \omega_\ell=\frac{\overline{P_\ell(x_\ell)}}{|P_\ell(x_\ell)|}\in\T_\K,
\end{equation}
where conjugation means the identity over $\R$, so that $\omega_\ell=\pm1$ in the real case. The denominator is nonzero, since $|P_\ell(x_\ell)|\ge M_\ell-\alpha_\ell^2\ge1-\alpha_\ell^2>0$. By this choice $\omega_\ell P_\ell(x_\ell)=|P_\ell(x_\ell)|$, so $d_\ell^{\omega_\ell}(x_\ell)=M_\ell-|P_\ell(x_\ell)|\le\alpha_\ell^2$ by \eqref{eq:nearmax}. Keeping this orientation fixed and applying \eqref{eq:backward} with $y=x_\ell$ at the stages $k,k+1,\ldots,\ell-1$, we obtain, for $k\le\ell$,
\[
 d_k^{\omega_\ell}(x_\ell)\le d_{k+1}^{\omega_\ell}(x_\ell)+\lambda_k\le\cdots\le d_\ell^{\omega_\ell}(x_\ell)+\sum_{h=k}^{\ell-1}\lambda_h ,
\]
where the sum is empty if $k=\ell$. Hence, by \eqref{eq:loss} and $\alpha_\ell\le\alpha_k$,
\begin{equation}\label{eq:base}
 d_k^{\omega_\ell}(x_\ell)
 \le\alpha_\ell^2+\sum_{h=k}^{\ell-1}\lambda_h
 \le5\alpha_k^2\qquad(k\le\ell).
\end{equation}

We claim that, for every $d\ge1$,
\begin{equation}\label{eq:compositions}
 \begin{split}
 &d_{i_1}^{\omega_\ell}
 \bigl(S_{i_1}(t_1)\cdots S_{i_d}(t_d)x_\ell\bigr)
 \le(7K_n)^d\alpha_{i_1}\\
 &\hspace{14mm}\text{for all } i_1<\cdots<i_d<\ell\text{ and all } t_1,\ldots,t_d\in[0,1].
 \end{split}
\end{equation}
As in Lemma~\ref{lem:geometry}, the operators act from right to left, and all the vectors involved lie in $B_X$ because every $S_k(t)$ is a contraction. By the composition formula for two averages given after \eqref{eq:Rfinite}, applying $\Rop_{i_1},\ldots,\Rop_{i_d}$ successively to a polynomial averages it over the points $S_{i_1}(t_1)\cdots S_{i_d}(t_d)x$.

Fix $\ell$. We prove the claim by induction on $d$, keeping the orientation $\omega_\ell$ fixed throughout the induction.

For the case $d=1$, let $i<\ell$ and $t\in[0,1]$. Since $i+1\le\ell$, \eqref{eq:base} applies at stage $i+1$ and gives $d_{i+1}^{\omega_\ell}(x_\ell)\le5\alpha_{i+1}^2=5\alpha_i^4\le5\alpha_i^2$. By \eqref{eq:insert} with $k=i$ and $y=x_\ell$, and by \eqref{eq:loss},
\[
 d_i^{\omega_\ell}(S_i(t)x_\ell)
 \le\frac{K_n}{\alpha_i}\bigl(d_{i+1}^{\omega_\ell}(x_\ell)+\lambda_i\bigr)
 \le\frac{K_n}{\alpha_i}(5\alpha_i^2+2\alpha_i^2)
 =7K_n\alpha_i.
\]

Now we prove the inductive step. Let $d\ge2$ and assume the claim for $d-1$. Let $i_1<\cdots<i_d<\ell$ and $t_1,\ldots,t_d\in[0,1]$, put $j=i_1$, $k=i_2$, $C_d=(7K_n)^{d-1}$, and
\[
 y=S_{i_2}(t_2)\cdots S_{i_d}(t_d)x_\ell\in B_X .
\]
The chain $i_2<\cdots<i_d<\ell$ has length $d-1$ and ends at the same vector $x_\ell$, so the induction hypothesis, applied with the parameters $t_2,\ldots,t_d$, gives $d_k^{\omega_\ell}(y)\le C_d\alpha_k$. Since $k\ge j+1$ and $(\alpha_h)$ is decreasing, $\alpha_k\le\alpha_{j+1}=\alpha_j^2$. Applying \eqref{eq:backward} at the fixed point $y$ at the stages $j+1,\ldots,k-1$, and then \eqref{eq:loss}, we get
\[
 d_{j+1}^{\omega_\ell}(y)
 \le d_k^{\omega_\ell}(y)+\sum_{h=j+1}^{k-1}\lambda_h
 \le C_d\alpha_k+4\alpha_{j+1}^2
 \le C_d\alpha_j^2+4\alpha_{j+1}^2
 =C_d\alpha_j^2+4\alpha_j^4 .
\]
The decay $\alpha_{j+1}=\alpha_j^2$ is used precisely here.
Now \eqref{eq:insert} at stage $j$, applied to this $y$, gives
\[
 \begin{split}
 d_j^{\omega_\ell}(S_j(t_1)y)
 &\le\frac{K_n}{\alpha_j}\bigl(d_{j+1}^{\omega_\ell}(y)+\lambda_j\bigr)
 \le\frac{K_n}{\alpha_j}\bigl(C_d\alpha_j^2+4\alpha_j^4+2\alpha_j^2\bigr)
 =K_n\bigl((C_d+2)\alpha_j+4\alpha_j^3\bigr)\\
 &\le K_n(C_d+6)\alpha_j\le7K_nC_d\alpha_j=(7K_n)^d\alpha_{i_1},
 \end{split}
\]
where we used $\alpha_j^3\le\alpha_j$ and $C_d\ge1$. Since $S_j(t_1)y=S_{i_1}(t_1)\cdots S_{i_d}(t_d)x_\ell$, and since the indices and the parameters $t_1,\ldots,t_d\in[0,1]$ were arbitrary, this proves \eqref{eq:compositions} for $d$.

Now let $i_1<\cdots<i_n$, put $j=i_1$, $\ell=i_n$, and write $\Gamma_n=(7K_n)^{n-1}$. Applying \eqref{eq:compositions} with $d=n-1$ to the chain $i_1<\cdots<i_{n-1}<\ell$, and recalling from \eqref{eq:deficits} that $d_j^{\omega_\ell}=M_j-\Rea(\omega_\ell P_j)$ on $B_X$, we obtain
\[
 \begin{split}
 M_j-\Rea\bigl(\omega_\ell P_j(S_{i_1}(t_1)\cdots S_{i_{n-1}}(t_{n-1})x_\ell)\bigr)
 &=d_j^{\omega_\ell}\bigl(S_{i_1}(t_1)\cdots S_{i_{n-1}}(t_{n-1})x_\ell\bigr)\\
 &\le(7K_n)^{n-1}\alpha_j=\Gamma_n\alpha_j ,
 \end{split}
\]
that is,
\[
 \Rea\bigl(\omega_\ell P_j(
 S_{i_1}(t_1)\cdots S_{i_{n-1}}(t_{n-1})x_\ell)\bigr)
 \ge M_j-\Gamma_n\alpha_j
 \quad\text{for every } t_1,\ldots,t_{n-1}\in[0,1].
\]
This is hypothesis \eqref{eq:path} of Lemma~\ref{lem:geometry} for the polynomial $\omega_\ell P_j$ and $M=M_j$ (note that $\norm{\omega_\ell P_j}=M_j\ge1$ because $|\omega_\ell|=1$), with the points $u_r=x_{i_r}$ $(1\le r\le n)$ and the functionals $f_{i_r}$ in place of $f_r$ $(1\le r<n)$, so that the maps $S_r(t)$ of the lemma are the maps $S_{i_r}(t)$, and with $E=\Gamma_n\alpha_j$. The symmetric form associated with $\omega_\ell P_j$ is $\omega_\ell A_j$, so Lemma~\ref{lem:geometry} gives $|\omega_\ell A_j(x_{i_1},\ldots,x_{i_n})|\ge M_j-H_n\sqrt{M_j\Gamma_n\alpha_j}$, that is,
\begin{equation}\label{eq:stageMixed}
 |A_j(x_{i_1},\ldots,x_{i_n})|
 \ge M_j-H_n\sqrt{M_j\Gamma_n\alpha_j}.
\end{equation}
In particular the phase $\omega_\ell$ has disappeared, and \eqref{eq:stageMixed} holds for every increasing $n$-tuple. Let $B=\check Q$, and recall that $M=\norm Q$. We showed after \eqref{eq:tailQ} that $\norm{B-A_j}\le4c_n\alpha_j$, and $|M-M_j|\le\norm{Q-P_j}\le4\alpha_j$. Since $x_{i_1},\ldots,x_{i_n}\in B_X$ and $M_j<2$,
\[
 \begin{split}
 |B(x_{i_1},\ldots,x_{i_n})|
 &\ge|A_j(x_{i_1},\ldots,x_{i_n})|-\norm{B-A_j}\\
 &\ge M_j-H_n\sqrt{M_j\Gamma_n\alpha_j}-4c_n\alpha_j\\
 &\ge M-4\alpha_j-H_n\sqrt{2\Gamma_n\alpha_j}-4c_n\alpha_j .
 \end{split}
\]
Thus
\begin{equation}\label{eq:ordered}
 |B(x_{i_1},\ldots,x_{i_n})|\ge M-\eta_{i_1}
 \qquad(i_1<\cdots<i_n),
\end{equation}
where
\begin{equation}\label{eq:eta}
 \eta_j=H_n\sqrt{2\Gamma_n\alpha_j}+4(c_n+1)\alpha_j
 \longrightarrow0.
\end{equation}
Here $\eta_j\to0$ because $\alpha_j=a^{2^{j-1}}\to0$ as $j\to\infty$, since $0<a<1$. The estimate \eqref{eq:ordered} holds for every strictly increasing $n$-tuple, and its error depends only on the first index $i_1$.

Choose a free ultrafilter $\U$ on $\N$, that is, an ultrafilter containing every cofinite subset of $\N$ (such an ultrafilter exists by Zorn's lemma). Every bounded scalar sequence $(a_k)$ has a unique limit $\ul ka_k$ along $\U$. This limit is linear in the sequence, satisfies $|\ul ka_k|=\ul k|a_k|$, preserves inequalities between real sequences which hold on a set belonging to $\U$, and coincides with $\lim_ka_k$ when the latter exists. Define
\begin{equation}\label{eq:z}
 z(f)=\ul{k}f(x_k)\qquad(f\in X^*),
\end{equation}
where $(x_k)$ is the sequence chosen in \eqref{eq:nearmax}. Since $|f(x_k)|\le\norm f$, the limit exists, $z$ is $\K$-linear in $f$, and $|z(f)|\le\norm f$; thus $z\in B_{X^{**}}$. Put $B_n=B$ and define $B_{r-1}$ from $B_r$ by \eqref{eq:partial}, using this $z$. For fixed $v_1,\ldots,v_{r-1}$, the functional $v\mapsto B_r(v_1,\ldots,v_{r-1},v)$ belongs to $X^*$, so \eqref{eq:z} gives
\begin{equation}\label{eq:ulpartial}
 B_{r-1}(v_1,\ldots,v_{r-1})
 =\ul{k}B_r(v_1,\ldots,v_{r-1},x_k),
 \qquad B_0=\widetilde Q(z),
\end{equation}
where the last equality is the identity $A_0=\widetilde P(z)$ established after \eqref{eq:partial}. Thus $\widetilde Q(z)$ is the iterated limit $\ul{i_1}\cdots\ul{i_n}B(x_{i_1},\ldots,x_{i_n})$, in which the limit in the last index is taken first.

We prove by descending induction on $r$ that
\[
 |B_r(x_{i_1},\ldots,x_{i_r})|\ge M-\eta_{i_1}
 \quad\text{for all } i_1<\cdots<i_r,\qquad r=n,\ldots,1.
\]
The case $r=n$ is \eqref{eq:ordered}. Let $2\le r\le n$, assume the estimate for $r$, and fix $i_1<\cdots<i_{r-1}$. For every $k>i_{r-1}$ the tuple $(i_1,\ldots,i_{r-1},k)$ is strictly increasing, so $|B_r(x_{i_1},\ldots,x_{i_{r-1}},x_k)|\ge M-\eta_{i_1}$. The set $\{k:k>i_{r-1}\}$ is cofinite and therefore belongs to $\U$. By \eqref{eq:ulpartial} and the properties of the limit along $\U$,
\[
 |B_{r-1}(x_{i_1},\ldots,x_{i_{r-1}})|
 =\ul k|B_r(x_{i_1},\ldots,x_{i_{r-1}},x_k)|\ge M-\eta_{i_1},
\]
which is the estimate for $r-1$. Finally, using the case $r=1$ and $\eta_k\to0$,
\[
 |\widetilde Q(z)|=|B_0|
 =\ul{k}|B_1(x_k)|\ge\ul{k}(M-\eta_k)=M.
\]
On the other hand, $z\in B_{X^{**}}$ and Lemma~\ref{lem:isometry} give $|\widetilde Q(z)|\le\norm{\widetilde Q}=\norm Q=M$. Hence $|\widetilde Q(z)|=\norm{\widetilde Q}=\norm Q$, and \eqref{eq:tailQ} gives $\norm{Q-P}<\eps$.
\end{proof}

The proof just given uses the same estimate over the two scalar fields. In the real case that estimate can be improved as follows.

\begin{remark}\label{rem:real}
Under the hypotheses of Lemma~\ref{lem:geometry}, if $\K=\R$, then
\begin{equation}\label{eq:realbound}
 |\check P(u_1,\ldots,u_n)|\ge M-C_nE,
 \qquad C_n=D_n\bigl(1+2n(n-1)c_n\bigr).
\end{equation}
For $E\ge M$ the assertion is trivial because $C_n\ge1$. For $0\le E<M$, retain the notation $\rho$, $\xi$, $\upsilon$, $y_i$ and $v_i$ from the proof of Lemma~\ref{lem:geometry}. For real numbers $\xi,\upsilon$ one has
\[
 |\xi-\upsilon|+|\xi+\upsilon|=2\max\{|\xi|,|\upsilon|\}\le2 :
\]
if $\xi$ and $\upsilon$ have the same sign, then $|\xi+\upsilon|=|\xi|+|\upsilon|$ and $|\xi-\upsilon|=\bigl||\xi|-|\upsilon|\bigr|$, and otherwise the two expressions exchange their roles. Since $|\xi+\upsilon|\ge2(1-\rho)^{1/n}$ by the first part of that proof, and $(1-\rho)^{1/n}\ge1-\rho$, it follows that
\[
 |\xi-\upsilon|\le2\bigl[1-(1-\rho)^{1/n}\bigr]\le2\rho.
\]
Thus \eqref{eq:compatibility} holds with $2\rho$ in place of $2\sqrt\rho$, and the comparison argument leading to \eqref{eq:convexdistance} gives $\norm{y_1-v_1}\le2(n-1)\rho$. Moreover $P(y_1)$ is real and $M-E\le P(y_1)\le M$, so $|P(y_1)-M|\le E$. Hence, by \eqref{eq:lipschitz} and $M\rho=E$,
\[
 |P(v)-M|\le nc_nM\cdot2(n-1)\rho+E=\bigl(1+2n(n-1)c_n\bigr)E
 \quad\bigl(v\in\conv\{z_1,\ldots,z_n\}\bigr)
\]
in place of \eqref{eq:flat}, and the finite-difference argument at the end of the same proof yields \eqref{eq:realbound}. In the real version of the main proof, \eqref{eq:stageMixed} is then replaced by $|A_j(x_{i_1},\ldots,x_{i_n})|\ge M_j-C_n\Gamma_n\alpha_j$, and one may take
\[
 \eta_j=\bigl(C_n\Gamma_n+4c_n+4\bigr)\alpha_j
\]
instead of \eqref{eq:eta}. The rest of the proof remains the same.
\end{remark}

\section{Vector-valued polynomials}\label{sec:vector}

Let $Y$ be a Banach space over $\K$. A map $P:X\to Y$ is a continuous $n$-homogeneous polynomial if $P(x)=L(x,\ldots,x)$ for some continuous $n$-linear map $L:X\times\cdots\times X\to Y$. We denote by $\PnY$ the space of all such polynomials, with the norm $\norm P=\sup_{x\in B_X}\norm{P(x)}$, and by $\check P$ the continuous symmetric $n$-linear map associated with $P$. For every $y^*\in Y^*$, the scalar polynomial $y^*\circ P$ belongs to $\Pn$, and its associated symmetric form is $y^*\circ\check P$. Since $\norm v=\sup_{y^*\in B_{Y^*}}|y^*(v)|$ for $v\in Y$, applying \eqref{eq:polarnorm} to the polynomials $y^*\circ P$ we obtain $\norm P\le\norm{\check P}\le c_n\norm P$. Moreover, $\PnY$ is a Banach space; see \cite{Di}.

The canonical (Aron--Berner) extension of $P\in\PnY$ is the map $\widetilde P:X^{**}\to Y^{**}$ defined by
\begin{equation}\label{eq:ABvector}
 \widetilde P(z)(y^*)=\widetilde{y^*\circ P}(z)\qquad(z\in X^{**},\ y^*\in Y^*).
\end{equation}
By Lemma~\ref{lem:isometry}, the map
$y^*\mapsto\widetilde{y^*\circ P}(z)$ is linear and satisfies
\[
|\widetilde{y^*\circ P}(z)|
\le \norm{y^*}\norm P\norm z^n.
\]
Thus $\widetilde P(z)\in Y^{**}$ and
$\norm{\widetilde P(z)}\le\norm P\norm z^n$.
To see that $\widetilde P$ is an $n$-homogeneous polynomial,
we construct an $n$-linear map
$\mathcal A:(X^{**})^n\to Y^{**}$.
For each $y^*\in Y^*$, let
$A_{y^*}=y^*\circ\check P$, and define
$\mathcal A(z_1,\ldots,z_n)\in Y^{**}$ by
\[
\mathcal A(z_1,\ldots,z_n)(y^*)
=\overline{A_{y^*}}(z_1,\ldots,z_n),
\qquad y^*\in Y^*.
\]
This definition is meaningful because the right-hand side
depends linearly on $y^*$ and
\[
|\overline{A_{y^*}}(z_1,\ldots,z_n)|
\le \norm{y^*}\norm{\check P}
\prod_{j=1}^n\norm{z_j}.
\]
Moreover, the $n$-linearity of each $\overline{A_{y^*}}$ implies that $\mathcal A$ is a continuous $n$-linear map. By the definition of the canonical extension, $\widetilde P(z)=\mathcal A(z,\ldots,z)$, so $\widetilde P$ is a continuous $n$-homogeneous polynomial. It extends $P$ when $Y$ is identified with its canonical image in $Y^{**}$, whence $\norm{\widetilde P}=\norm P$. The map $P\mapsto\widetilde P$ is linear, and for $n=1$ it coincides with the second adjoint $P^{**}$.

\begin{theorem}\label{thm:vector}
Let $X$ and $Y$ be Banach spaces over $\K$ and let $n\ge1$. The set of all $P\in\PnY$ whose canonical extensions attain their norms is dense in $\PnY$. More precisely, for every $P\in\PnY$ and every $\eps>0$ there exist $Q\in\PnY$ and $z\in B_{X^{**}}$ such that
\[
 \norm{Q-P}<\eps,
 \qquad
 \norm{\widetilde Q(z)}=\norm{\widetilde Q}=\norm Q.
\]
\end{theorem}

For $Y=\K$ this is Theorem~\ref{thm:main}, and for $n=1$ it is Lindenstrauss' theorem \cite[Theorem~1]{L}.

The proof of Theorem~\ref{thm:main} does not extend to Theorem~\ref{thm:vector} simply by replacing absolute values by norms. In the last step of that proof, the lower bounds for $|\check Q(x_{i_1},\ldots,x_{i_n})|$ survive the iterated limits along $\U$ because the modulus is continuous. For $Y$-valued polynomials, $\widetilde Q(z)$ is obtained from the vectors $\check Q(x_{i_1},\ldots,x_{i_n})$ through iterated weak-star limits in $Y^{**}$, and the norm is only weak-star lower semicontinuous. What is needed are lower bounds for $|y^*(\check Q(x_{i_1},\ldots,x_{i_n}))|$ with a \emph{fixed} functional $y^*\in B_{Y^*}$; in the proof of Theorem~\ref{thm:main}, however, the orientation $\omega_\ell$ depends on the largest index $\ell=i_n$, which is the one that varies in the innermost limit. We therefore modify the perturbations so that norming functionals of $Y$ enter the construction; the induction then transfers the functional to the smallest index of each chain.

\begin{proof}[Proof of Theorem~\ref{thm:vector}]
For $n=1$ the theorem is \cite[Theorem~1]{L}. The case $P=0$ is trivial and, exactly as in the proof of Theorem~\ref{thm:main}, the linearity of $P\mapsto\widetilde P$ allows us to assume that $\norm P=1$. So let $n\ge2$, $\norm P=1$ and $\eps>0$, and choose $a$ and $\alpha_k$ as in \eqref{eq:alphas}; then \eqref{eq:tails} holds.

Starting with $P_1=P$, we construct inductively $P_k\in\PnY$, $x_k\in S_X$, $f_k\in S_{X^*}$ and $y_k^*\in S_{Y^*}$ such that, with $M_k=\norm{P_k}$,
\[
 \norm{P_k(x_k)}\ge M_k-\alpha_k^2,
 \qquad f_k(x_k)=1,
 \qquad y_k^*(P_k(x_k))=\norm{P_k(x_k)},
\]
and
\begin{equation}\label{eq:viteration}
 P_{k+1}=P_k+\alpha_k(\Rop_kp_k)\,w_k,
 \qquad\text{where}\quad
 p_k=y_k^*\circ P_k\in\Pn,
 \quad
 w_k=\frac{P_k(x_k)}{\norm{P_k(x_k)}},
\end{equation}
with $S_k(t)=S_{x_k,f_k}(t)$ and $\Rop_k=\Rop_{x_k,f_k}$ as before. The vectors $x_k$ and the functionals $f_k$ exist as in the proof of Theorem~\ref{thm:main}, and $y_k^*$ exists by the Hahn--Banach theorem. The vector $w_k$ is well defined because $M_k\ge1$ (this follows inductively from \eqref{eq:normbounds}, as in the proof of Theorem~\ref{thm:main}), so that $\norm{P_k(x_k)}\ge1-\alpha_k^2>0$. Note that $\norm{p_k}\le M_k$, $\norm{w_k}=1$ and $\Rop_kp_k(x_k)=p_k(x_k)=\norm{P_k(x_k)}$. Hence, by \eqref{eq:Rproperties},
\[
 \norm{P_{k+1}-P_k}=\alpha_k\norm{\Rop_kp_k}\le\alpha_kM_k,
 \qquad
 P_{k+1}(x_k)=(1+\alpha_k)\norm{P_k(x_k)}\,w_k .
\]
Therefore \eqref{eq:increment} and \eqref{eq:normbounds} hold, and the arguments following them apply without change: $M_k\ge1$ for every $k$, \eqref{eq:product} holds, $(P_k)$ converges in $\PnY$ to some $Q$ satisfying \eqref{eq:tailQ}, $M=\norm Q$ satisfies $1\le M<2$, $\norm{\check Q-\check P_k}\le4c_n\alpha_k$, and the numbers $\lambda_k=(1+\alpha_k)M_k-M_{k+1}$ satisfy \eqref{eq:loss}.

For $g\in B_{Y^*}$ and $y\in B_X$ put
\[
 d_k^g(y)=M_k-\Rea g(P_k(y)).
\]
Then $d_k^g\ge0$ on $B_X$, because $|g(P_k(y))|\le\norm{P_k(y)}\le M_k$. Let $g\in B_{Y^*}$ and $y\in B_X$. By \eqref{eq:viteration},
\[
 g(P_{k+1}(y))=g(P_k(y))+\alpha_k\,g(w_k)\int_0^1p_k(S_k(t)y)\,dt,
\]
and $g(w_k)p_k=h\circ P_k$ with $h=g(w_k)y_k^*\in B_{Y^*}$. The computation in \eqref{eq:deficitidentity} therefore gives
\begin{equation}\label{eq:vdeficit}
 d_k^g(y)+\alpha_k\int_0^1d_k^{\,g(w_k)y_k^*}(S_k(t)y)\,dt=d_{k+1}^g(y)+\lambda_k .
\end{equation}
Both terms on the left-hand side are nonnegative, and $t\mapsto d_k^{\,g(w_k)y_k^*}(S_k(t)y)$ is a real polynomial of degree at most $n$, nonnegative on $[0,1]$. As in the proof of \eqref{eq:backward} and \eqref{eq:insert}, we obtain
\begin{align}
 d_k^g(y)&\le d_{k+1}^g(y)+\lambda_k,\label{eq:vbackward}\\
 \sup_{t\in[0,1]}d_k^{\,g(w_k)y_k^*}(S_k(t)y)
 &\le\frac{K_n}{\alpha_k}\bigl(d_{k+1}^g(y)+\lambda_k\bigr).\label{eq:vinsert}
\end{align}
Thus, in contrast with \eqref{eq:insert}, the insertion of the substitution $S_k$ replaces the functional $g$ by $g(w_k)y_k^*$, a multiple of the functional $y_k^*$ chosen at stage $k$.

Fix $\ell\in\N$. Since $y_\ell^*(P_\ell(x_\ell))=\norm{P_\ell(x_\ell)}\ge M_\ell-\alpha_\ell^2$, we have $d_\ell^{\,y_\ell^*}(x_\ell)\le\alpha_\ell^2$, and \eqref{eq:vbackward} with $g=y_\ell^*$ gives, exactly as in the proof of \eqref{eq:base},
\begin{equation}\label{eq:vbase}
 d_k^{\,y_\ell^*}(x_\ell)\le5\alpha_k^2\qquad(k\le\ell).
\end{equation}
We claim that, for every $d\ge1$ and all indices $i_1<\cdots<i_d<\ell$, there is a scalar $c$ with $|c|\le1$, depending only on $i_1,\ldots,i_d$ and $\ell$, such that
\begin{equation}\label{eq:vcompositions}
 d_{i_1}^{\,c\,y_{i_1}^*}\bigl(S_{i_1}(t_1)\cdots S_{i_d}(t_d)x_\ell\bigr)\le(7K_n)^d\alpha_{i_1}
 \quad\text{for all } t_1,\ldots,t_d\in[0,1].
\end{equation}
For $d=1$ and $i<\ell$, apply \eqref{eq:vinsert} with $k=i$, $g=y_\ell^*$ and $y=x_\ell$; since $d_{i+1}^{\,y_\ell^*}(x_\ell)\le5\alpha_{i+1}^2\le5\alpha_i^2$ by \eqref{eq:vbase} and $\lambda_i\le2\alpha_i^2$, we obtain \eqref{eq:vcompositions} with $c=y_\ell^*(w_i)$. Now let $d\ge2$ and assume the claim for $d-1$. Let $i_1<\cdots<i_d<\ell$, put $j=i_1$, $k=i_2$ and $C_d=(7K_n)^{d-1}$, fix $t_2,\ldots,t_d\in[0,1]$, and let $y=S_{i_2}(t_2)\cdots S_{i_d}(t_d)x_\ell\in B_X$. By the induction hypothesis there is a scalar $c'$ with $|c'|\le1$, independent of $t_2,\ldots,t_d$, such that $d_k^{\,g}(y)\le C_d\alpha_k$ for $g=c'y_k^*\in B_{Y^*}$. Applying \eqref{eq:vbackward} with this fixed $g$ at the stages $j+1,\ldots,k-1$, and then \eqref{eq:loss}, we obtain $d_{j+1}^{\,g}(y)\le C_d\alpha_j^2+4\alpha_j^4$, exactly as in the proof of \eqref{eq:compositions}. Then \eqref{eq:vinsert} at stage $j$ gives, for every $t_1\in[0,1]$,
\[
 d_j^{\,g(w_j)y_j^*}(S_j(t_1)y)\le\frac{K_n}{\alpha_j}\bigl(C_d\alpha_j^2+4\alpha_j^4+2\alpha_j^2\bigr)\le7K_nC_d\alpha_j=(7K_n)^d\alpha_{i_1}.
\]
Since $S_j(t_1)y=S_{i_1}(t_1)\cdots S_{i_d}(t_d)x_\ell$, this is \eqref{eq:vcompositions} with $c=g(w_j)=c'y_k^*(w_j)$, which satisfies $|c|\le1$ and does not depend on $t_1,\ldots,t_d$. This proves the claim, for every $\ell\in\N$.

Now let $m<i_1<\cdots<i_n$ and put $\ell=i_n$. We apply \eqref{eq:vcompositions} with $d=n$ to the chain $m<i_1<\cdots<i_{n-1}<\ell$, and we take the parameter of $S_m$ equal to $0$. Since $S_m(0)$ is the identity and $d_m^{\,c\,y_m^*}=M_m-\Rea(c\,p_m)$ on $B_X$, we obtain a scalar $c$ with $|c|\le1$ such that
\[
 \Rea\bigl(c\,p_m(S_{i_1}(t_1)\cdots S_{i_{n-1}}(t_{n-1})x_\ell)\bigr)\ge M_m-(7K_n)^n\alpha_m
 \quad\text{for every } t_1,\ldots,t_{n-1}\in[0,1].
\]
The scalar polynomial $c\,p_m$ satisfies $\norm{c\,p_m}\le M_m$, and its associated symmetric form is $c\,y_m^*\circ\check P_m$. Hence Lemma~\ref{lem:geometry}, applied to $c\,p_m$ with $M=M_m$, the points $u_r=x_{i_r}$ $(1\le r\le n)$, the functionals $f_{i_r}$ $(1\le r<n)$, so that the maps $S_r(t)$ of the lemma are the maps $S_{i_r}(t)$, and $E=(7K_n)^n\alpha_m$, gives
\[
 |y_m^*(\check P_m(x_{i_1},\ldots,x_{i_n}))|
 \ge|c\,y_m^*(\check P_m(x_{i_1},\ldots,x_{i_n}))|
 \ge M_m-H_n\sqrt{M_m(7K_n)^n\alpha_m}.
\]
Let $B=\check Q$. Since $\norm{y_m^*}=1$, $\norm{B-\check P_m}\le4c_n\alpha_m$, $M_m<2$ and $M_m\ge M-4\alpha_m$, we conclude that
\begin{equation}\label{eq:vordered}
 |y_m^*(B(x_{i_1},\ldots,x_{i_n}))|\ge M-\eta'_m
 \qquad(m<i_1<\cdots<i_n),
\end{equation}
where $\eta'_m=H_n\sqrt{2(7K_n)^n\alpha_m}+4(c_n+1)\alpha_m$, so that $\eta'_m\to0$ as $m\to\infty$. In contrast with \eqref{eq:ordered}, the functional $y_m^*$ is now fixed, and its index is smaller than the indices of all the points.

Let $\U$ be a free ultrafilter on $\N$ and define $z\in B_{X^{**}}$ by \eqref{eq:z}, now with the vectors $x_k$ chosen above. Fix $m$, and let $b=y_m^*\circ B$, which is the symmetric form associated with the scalar polynomial $y_m^*\circ Q$. Put $b_n=b$ and define $b_{r-1}$ from $b_r$ by \eqref{eq:partial}, using this $z$. Descending induction on $r$, starting from \eqref{eq:vordered} and arguing exactly as after \eqref{eq:ulpartial} but with all the indices larger than $m$, gives $|b_r(x_{i_1},\ldots,x_{i_r})|\ge M-\eta'_m$ for all $m<i_1<\cdots<i_r$ and $r=n,\ldots,1$. Since $b_0=\widetilde{y_m^*\circ Q}(z)$ and $\{k:k>m\}\in\U$, we get
\[
 |\widetilde Q(z)(y_m^*)|=|\widetilde{y_m^*\circ Q}(z)|=|b_0|=\ul k|b_1(x_k)|\ge M-\eta'_m ,
\]
where the first equality is \eqref{eq:ABvector}. Hence $\norm{\widetilde Q(z)}\ge M-\eta'_m$ for every $m$, and therefore $\norm{\widetilde Q(z)}\ge M$. Since $z\in B_{X^{**}}$ and $\norm{\widetilde Q}=\norm Q=M$, we conclude that $\norm{\widetilde Q(z)}=\norm{\widetilde Q}=\norm Q$, while $\norm{Q-P}<\eps$ by \eqref{eq:tailQ}.
\end{proof}

When $Y=\K$, the functional $y_k^*$ is multiplication by $\overline{w_k}$, where $w_k=P_k(x_k)/|P_k(x_k)|$ is a unimodular scalar; hence $(\Rop_kp_k)w_k=\Rop_kP_k$, and \eqref{eq:viteration} reduces to \eqref{eq:iteration}. Thus the proof above also gives a variant of the last part of the proof of Theorem~\ref{thm:main}, in which the orientation is attached to the smallest index of each chain.

\medskip
\noindent\textbf{Acknowledgment.} The author acknowledges the use of artificial intelligence--based tools for language editing and stylistic improvements of the manuscript. The author remains fully responsible for the content of this work.

\begin{samepage}

\end{samepage}
\end{document}